\documentclass[12pt]{article}
\usepackage{amscd,amsfonts,amssymb,amsmath,latexsym,array,hhline,amsthm}
\newcommand{\paragr}{\hspace{6mm}}
\newtheorem{Theorem}{Theorem}[section]

\theoremstyle{remark}

\begin{document}

\sloppy

\author{Pavel G.~Grigoriev and Stanislav A. Molchanov\\
\small
Department of Mathematics and Statistics\\
\small
University of North Carolina at Charlotte\\
\small
Charlotte, NC 28223, USA
\\
\small
thepavel@mail.ru}

\title{On Decoupling of Functions of Normal Vectors II}

\date{}

\markboth{Pavel G.~Grigoriev, Stanislav Molchanov}{On Decoupling of Functions of Normal Vectors}

\maketitle

\begin{abstract}
{\small
A decoupling type inequality
for a sum of functions of  Guassian vectors is established.
}
\end{abstract}


\section{Results}

In \cite{grigmol1}
the following decoupling results were established.


\begin{Theorem}\label{TheoremA}
For a normally distributed random vector $\bar Y=(Y_i)_{i=1,\dots,d}$
satisfying $\mathsf{E} Y_i=0$, $\mathsf{E} Y_i^2=1$, $i=1,\dots,d$,
we have
\begin{equation*}\label{theoremold}
c^-\sum_{i=1}^d \|\varphi_i(Y_i)\|_2^2\le
\Big\|\sum_{i=1}^d \varphi_i(Y_i)\Big\|_2^2\le
c^+\sum_{i=1}^d \|\varphi_i(Y_i)\|_2^2
\end{equation*}
for all measurable functions $\varphi_i:\mathbb{R}\to\mathbb{R}$
satisfying $\mathsf{E} \varphi_i(Y_i)=0$
with the constants
$c^-$ and $c^+$ being
the smallest and the largest eigenvalues of the correlation matrix of $\bar Y$.
Moreover the constants are the best possible.
\end{Theorem}

\begin{Theorem}\label{TheoremDecoup}
Let $\bar Y_1=(Y_{1,i})_{i=1,\dots,d_1}$
and $\bar Y_2=(Y_{2,j})_{j=1,\dots,d_2}$
be standard normal vectors with
the correlations
\begin{equation*}\label{correlation-2vectro}
\mathsf{E} Y_{\alpha,i}Y_{\beta,j}=
\left\{
\begin{array}{ll}
0, &\alpha= \beta, i\neq j\\
1, &\alpha=\beta, i=j\\
\rho_{i,j}, &\alpha\neq \beta
\end{array}
\right.
\end{equation*}
Then
\begin{equation*}\label{naive2}
c_-\big(\|\varphi_1(\bar Y_1)\|_2^2+\|\varphi_2(\bar Y_2)\|_2^2\big)
\le
\big\|\varphi_1(\bar Y_1)+\varphi_2(\bar Y_2)\big\|_2^2\le
c_+\big(\|\varphi_1(\bar Y_1)\|_2^2+\|\varphi_2(\bar Y_2)\|_2^2\big)
\end{equation*}
holds for all measurable functions $\varphi_\alpha:\mathbb{R}^{d_\alpha}\to\mathbb{R}$
 such that $\mathsf{E}\varphi_\alpha(\bar Y_\alpha)=0$, $\alpha=1,2$, with
the  constants
$
c_\pm =1\pm s^*,
$
where $s^*$ is the maximum singular value
of the matrix
$R=(\rho_{i,j})_{i=1,\dots,d_1,j=1,\dots,d_2}$.
These constants cannot be improved.
\end{Theorem}

These theorems refine the estimates used by Cherny et al. \cite{cdm1}, \cite{cdm2}.
Here we generalize these results further and prove

\begin{Theorem}\label{TheoremNew}
Let $\bar Y_\alpha=(Y_{\alpha,i})_{i=1,\dots,d_\alpha}$,
$\alpha=1,\dots,N$,
be
standard normal vectors with
the correlations
\begin{equation*}\label{correlation-2vectro}
\mathsf{E} Y_{\alpha,i}Y_{\beta,j}=
\left\{
\begin{array}{ll}
0, &\alpha= \beta, i\neq j\\
1, &\alpha=\beta, i=j\\
\rho_{i,j}^{\alpha,\beta}, &\alpha\neq \beta
\end{array}
\right.
\end{equation*}
Then
\begin{equation}\label{newdecoup}
C_-
\sum_{\alpha=1}^N
\|\varphi_\alpha(\bar Y_\alpha)\|_2^2
\le
\Big\|\sum_{\alpha=1}^N
\varphi_\alpha(\bar Y_\alpha)\Big\|_2^2\le
C_+
\sum_{\alpha=1}^N
\|\varphi_\alpha(\bar Y_\alpha)\|_2^2
\end{equation}
for all measurable functions $\varphi_\alpha:\mathbb{R}^{d_\alpha}\to\mathbb{R}$,
satisfying $\mathsf{E}\varphi_\alpha(\bar Y_\alpha)=0$
 for all $\alpha=1,\dots,N$.
The  constants
$
C_\pm=1\pm \sigma_0
$, where
$\sigma_0$ denotes
the largest eigenvalue
of the matrix
$S^*=(s^*_{\alpha,\beta})_{1\le \alpha,\beta\le N}$
with
$s^*_{\alpha,\beta}$
being the maximum singular values
of the matrices
$R^{\alpha,\beta}=(\rho_{i,j}^{\alpha,\beta})_{i=1,\dots,d_\alpha,j=1,\dots,d_\beta}$
for $\alpha\neq\beta$
and $s^*_{\alpha,\alpha}:=0$.
\end{Theorem}

Note that while Theorems~\ref{TheoremA} and~\ref{TheoremDecoup}
give the decoupling estimates with sharp constants,
Theorem~\ref{TheoremNew} provides just rough estimate
with not the best constants (in particular $C_-$
could be negative). However, in view
of the applications described in
Cherny et al. \cite{cdm1}, \cite{cdm2}
the upper bound in~(\ref{newdecoup})
is still interesting.

\section{Proof of Theorem \ref{TheoremNew}}

The proof of Theorem~\ref{TheoremNew} follows the framework used in~\cite{grigmol1}.
We need to introduce notations used in~\cite{grigmol1}
and borrowed from~\cite{malmin}.

For $k=0,1,\dots$ we define
Wick's polynomials $:x^k:$ by
the extension
$$
\exp\Big(ax-\frac a2\Big)=\sum_{k=0}^\infty a^k\frac{:x^k:}{k!}.
$$
(Wick's polynomials are specially normalized Hermite's polynomials
used in mathematical physics.
We find these notations convenient for multidimensional case.)

Let
$\bar k=(k_1,\dots,k_d)$ be a $d$-dimensional vector of non-negative integers. Set
\begin{align}
\notag
|\bar k| &:=k_1+\dots+k_d,\\
\notag
\bar k! &:=k_1! k_2!\dots k_d!\\
\notag
\bar a^{\bar k} &:=a_1^{k_1}a_2^{k_2}\dots a_d^{k_d},
\qquad\text{for\ }\bar a :=(a_i)_{i=1,\dots,d}\in\mathbb{R}^d.
\end{align}
For a vector variable $\bar x=(x_k)_{k=1}^d$
we define multidimensional Wick's polynomial
by
$$
:\bar x^{\bar k}:\,
:=\prod_{i=1}^d\,:x_{i}^{k_i}:.
$$

It is well known (see e.g. \cite{malmin})
that for
a standard $d$-dimensional normal vector $\bar Y$
the system
$\big\{\bar k!^{-1/2}:\bar Y^{\bar k}:\big\}_{\bar k\in \mathbb{Z}_0^{d}}$
is an orthonormal bases
in the $L_2$ space generated by all square-integrable $f(\bar Y)$.
So for each $\alpha=1,\dots,N$, we have
\begin{equation*}
\label{align-theodec1}
\varphi_\alpha(\bar Y_\alpha) =\sum_{\bar k\in \mathbb{Z}_0^{d_\alpha}}
a_{\alpha, \bar k}\frac{:\bar Y_\alpha^{\bar k}:}{{\bar k!}^{1/2}}
=
\sum_{n=0}^\infty
\sum_{\bar k\in \mathbb{Z}_0^{d_\alpha},\ |\bar k|=n}
a_{\alpha, \bar k}\frac{:\bar Y_\alpha^{\bar k}:}{{\bar k!}^{1/2}}
\end{equation*}
and therefore
$$
\sum_{\alpha=1}^N
\varphi_\alpha(\bar Y_\alpha)
=
\sum_{\alpha=1}^N
\sum_{n=0}^\infty
\sum_{\bar k_\alpha\in \mathbb{Z}_0^{d_\alpha},\ |\bar k_\alpha|=n}
a_{\alpha, \bar k_\alpha}\frac{:\bar Y_\alpha^{\bar k_\alpha}:}{{\bar k_\alpha!}^{1/2}}
=\sum_{n=0}^\infty\sum_{\alpha=1}^N P_n(\bar Y_\alpha),
$$
where
$P_n(\bar Y_\alpha):=
\sum\limits_{\bar k_\alpha\in \mathbb{Z}_0^{d_\alpha},\ |\bar k_\alpha|=n}
a_{\alpha, \bar k_\alpha}\frac{:\bar Y_\alpha^{\bar k_\alpha}:}{{\bar k_\alpha!}^{1/2}}
$.
In particular,
$$
\sum_{\alpha=1}^N
\|\varphi_\alpha(\bar Y_\alpha)\|_2^2
=
\sum_{n=0}^\infty\sum_{\alpha=1}^N \|P_n(\bar Y_\alpha)\|_2^2.
$$

It is well-known that
$P_{n_1}(\bar Y_\alpha)$
and
$P_{n_2}(\bar Y_\beta)$
are orthogonal whenever $n_1\neq n_2$ (see e.g. \cite{grigmol1}).
Consequently, to prove
(\ref{newdecoup}) for arbitrary $\varphi_\alpha$
it suffices to prove it for
$\phi_\alpha(\bar Y_\alpha)=P_{n}(\bar Y_\alpha)$ for each $n=0,1,\dots$.

For a fixed $n$, we have
\begin{equation}\label{Pn}
\Big\|\sum_{\alpha=1}^N P_{n}(\bar Y_\alpha)\Big\|_2^2
=
\sum_{\alpha=1}^N
\|P_{n}(\bar Y_\alpha)\|_2^2
+
2
\sum_{\alpha<\beta}
\mathsf{E}
\big[P_{n}(\bar Y_\alpha)P_{n}(\bar Y_\beta)\big].
%
%
%
%
\end{equation}

Fix a pair $\alpha<\beta$.
Let $R^{\alpha,\beta}=U\Sigma V^T$ be the singular
value decomposition
of the matrix $R^{\alpha,\beta}$
(recall that here $\Sigma$ is a $d_\alpha\times d_\beta$ diagonal
matrix whose diagonal entries are the singular values of $R^{\alpha,\beta}$ and $U$, $V$ are orthogonal matrices of
corresponding sizes).
Let
$\bar Z_1:=U \bar Y_\alpha$
and
$\bar Z_2:=V \bar Y_\beta$
(all the vectors are assumed being columns). This transformation can be written in the
block matrix form
$$
\left(
\begin{array}{c}
\bar Z_1\\
\bar Z_2
\end{array}
\right)
=
\left(
\begin{array}{cc}
U &0\\
0 &V
\end{array}
\right)
\left(
\begin{array}{c}
\bar Y_\alpha\\
\bar Y_\beta
\end{array}
\right)
$$
with the obviously {\it orthogonal} $(d_\alpha+d_\beta)\times(d_\alpha+d_\beta)$ transformation matrix.
It follows that
\begin{equation}\label{YtoZ}
\mathsf{E}
\big[P_{n}(\bar Y_\alpha)P_{n}(\bar Y_\beta)\big]
=
\mathsf{E}
\big[P_{n}(\bar Z_1)P_{n}(\bar Z_2)\big]
=
\sum_{\bar k_\alpha\in \mathbb{Z}_0^{d_\alpha}\atop |\bar k_\alpha|=n}
\sum_{\bar k_\beta\in \mathbb{Z}_0^{d_\beta}\atop |\bar k_\beta|=n}
a_{\alpha, \bar k_\alpha}
a_{\beta, \bar k_\beta}
\mathsf{E}
\frac{:\bar Z_1^{\bar k_\alpha}:}{{\bar k_\alpha!}^{1/2}}\frac{:\bar Z_2^{\bar k_\beta}:}{{\bar k_\beta!}^{1/2}}
\end{equation}

Note that the
covariance structure of
$(\bar Z_1, \bar Z_2)$
is relatively simple, its covariance matrix is
$$
\left(
\begin{array}{cc}
I_{d_\alpha} &\Sigma\\
\Sigma^T &I_{d_\beta}
\end{array}
\right).
$$

Without loss of generality we assume that $d_\alpha\ge d_\beta$, i.e. the vectors $\bar Y_\alpha$
are ordered according to their dimensions.
Let us agree that the index vectors of different
dimensions are equal ($\bar k_\alpha= \bar k_\beta$) if
the shorter vector $\bar k_\beta$ coincides  with the first $d_\beta$ entries
of $k_\alpha$ and the other entries of $k_\alpha$ are zeros.
Also let us denote by $\bar s$ the vector of the diagonal entries of $\Sigma$
(i.e. the singular values of $R^{\alpha,\beta}$).

Using Lemma~3.1 from~\cite{grigmol1}
(for the case of diagonal $R$, see also formula (3.5) in~\cite{grigmol1})
we conclude
\begin{equation}\label{covforZ}
\mathsf{E} :\bar Z_1^{\bar k_\alpha}::\bar Z_2^{\bar k_\beta}:
=
\left\{
\begin{array}{ll}
0, &\bar k_\alpha\neq \bar k_\beta\\
\bar k_\beta!
\bar s^{\bar k_\beta},
&\bar k_\alpha= \bar k_\beta
\end{array}
\right.
\end{equation}

Using (\ref{covforZ}) we continue
(\ref{YtoZ}) as follows
\begin{equation*}
\mathsf{E}
\big[P_{n}(\bar Y_\alpha)P_{n}(\bar Y_\beta)\big]
=
\sum_{\bar k_\beta\in \mathbb{Z}_0^{d_\beta}, \bar k_\alpha\in \mathbb{Z}_0^{d_\alpha}:
\atop |\bar k_\beta|=n,\,k_\alpha= \bar k_\beta}
a_{\alpha, \bar k_\alpha}
a_{\beta, \bar k_\beta}
\bar s^{\bar k_\beta}.
\end{equation*}
Taking into account that the singular values lie in $[0,1]$
(because they are certain correlations)
we can estimate
\begin{align*}
\big|\mathsf{E}
\big[P_{n}(\bar Y_\alpha)P_{n}(\bar Y_\beta)\big]\big|
&\le s_{\alpha,\beta}^*
\sum_{\bar k_\beta\in \mathbb{Z}_0^{d_\beta}, \bar k_\alpha\in \mathbb{Z}_0^{d_\alpha}:
\atop |\bar k_\beta|=n,\,k_\alpha= \bar k_\beta}
|a_{\alpha, \bar k_\alpha}
a_{\beta, \bar k_\beta}|\\
&\le
s_{\alpha,\beta}^*
\Big(
\sum_{\bar k_\alpha\in \mathbb{Z}_0^{d_\alpha}\atop |\bar k_\alpha|=n}
a_{\alpha, \bar k_\alpha}^2
\Big)^{\frac12}
\Big(
\sum_{\bar k_\beta\in \mathbb{Z}_0^{d_\beta}\atop |\bar k_\beta|=n}
a_{\beta, \bar k_\beta}^2
\Big)^{\frac12}
=
s_{\alpha,\beta}^*
\|P_n(\bar Y_\alpha)\|_2\|P_n(\bar Y_\beta)\|_2
\end{align*}

Recalling (\ref{Pn}) we conclude
$$
\Big|
\Big\|\sum_{\alpha=1}^N P_{n}(\bar Y_\alpha)\Big\|_2^2
-
\sum_{\alpha=1}^N
\|P_{n}(\bar Y_\alpha)\|_2^2
\Big|
\le
2
\sum_{\alpha<\beta}
s_{\alpha,\beta}^*
\|P_n(\bar Y_\alpha)\|_2\|P_n(\bar Y_\beta)\|_2.
$$
By standard linear algebra argument
we have that the right hand side is bounded by
$\sigma_0\sum_{\alpha=1}^N \|P_n(\bar Y_\alpha)\|_2^2$,
where $\sigma_0$ is the largest
eigenvalue of the matrix $S^*$.
We conclude
that for a fixed $n$ the estimate
(\ref{newdecoup})
holds
for
$\phi_\alpha(\bar Y_\alpha)=P_{n}(\bar Y_\alpha)$, $\alpha=1,\dots,N$,
with $C_\pm=1\pm \sigma_0$.
As it was already pointed out this means that
we have (\ref{newdecoup}) for all $\phi_\alpha$ with the same constants.
$\Box$

%

%


\end{document}